\documentclass[reqno, 12pt]{amsart}
\usepackage{amsmath}
\usepackage{latexsym}
\usepackage{amsfonts}
\usepackage{amssymb}
\usepackage{amsfonts,euscript}
\usepackage{graphicx}
\DeclareGraphicsExtensions{.eps,epsfig,.bmp,.jpg,.pdf,.mps,.png,.gif}

 \usepackage[pagewise]{lineno}

\theoremstyle{plain}
\newtheorem{theorem}{\bf Theorem}[section]

\newtheorem{lemma}{\bf Lemma}[section]
\theoremstyle{definition}

\theoremstyle{remark}
\newtheorem{remark}{\bf Remark}[section]

\numberwithin{equation}{section}

\title[Stabilization of 1-D  parabolic equations with nonlocal B.C.]{Exponential stabilization of  the semilinear heat equation with nonlocal boundary conditions}
\author{Ionu\c{t} Munteanu}
\address{\noindent Ionu\c{t} Munteanu\newline Alexandru Ioan Cuza University, Department of Mathematics, and Octav Mayer Institute of Mathematics (Romanian Academy) 
\newline  700506 Ia\c{s}i, Romania}
 \email{ionut.munteanu@uaic.ro}

\begin{document}
\maketitle
\begin{abstract} The present work is devoted to the problem of boundary stabilization of the semilinear 1-D heat equation with nonlocal boundary conditions. The stabilizing controller is finite-dimensional, linear, given in an explicit form, involving only the eigenfunctions of the Laplace operator with nonlocal boundary conditions. 
	\end{abstract}

\section{Introduction of the problem}
The subject of the present paper is the following semi-linear heat equation with integral boundary conditions
\begin{equation}\label{oe1} \left\lbrace \begin{array}{l}y_t(t,x)-y''(t,x)-cy(t,x)=f(t,x,y(t,x)),\\
 t>0, x\in(0,\pi),\\
\\
\int_0^\pi \left[ y(t,x)dh_1(x)+y'(t,x)dh_2(x)\ \right]=0,\ t\geq0, \\
\\
 \int_0^\pi\left[y(t,x)dh_3(x)+y'(t,x)dh_4(x) \right]=0,\ t\geq0,\\  
 \\
y(0,x)=y_o(x),\ x\in(0,\pi).
\end{array} \right. \end{equation}
Here, $h_i,\ i=1,2,3,4,$ are functions of bounded variations, and the integration is in the sense of Riemann-Stieltjes, and $c>0$. The nonlinear function $f$ satisfies
$$|f(t,x,y)|\leq C |y|^m,\ \forall t\geq0,\ x\in (0,\pi),\ y\in \mathbb{R},$$where $C$ is some positive constant and $m>0$. And also
$$f(t,x,0)=0.$$  

Boundary-value problems, with integral boundary conditions, constitute a
very interesting class of problems, because they include as special cases two,
three, multi-point and nonlocal boundary-value problems. They arise naturally in thermal conduction, semiconductor  or hydrodynamic problems. For details see \cite{bouch} and the references therein.

We are concerned here with the problem of exponential asymptotic stability of the null solution, namely whether the equation \eqref{oe1} is well-posed and the corresponding solution satisfies
\begin{equation}\label{s1}\int_0^\pi y^2(t,x)dx\leq C e^{-\rho t}\int_0^\pi y^2_o(x)dx,\ \forall t\geq0,\end{equation} for some positive constants $C,\rho.$   In \cite{c15}, it is shown that the finite time blow-up phenomenon may occur in the case of  semilinear heat equation with nonlocal boundary conditions. Besides this, the positive constant $c$ may have an important contribution to the instability. Therefore it makes sense to pose the problem of  stabilization of equations of the type \eqref{oe1}. That is, we shall see the functions $h_1,h_2$ as controllers. More precisely, given $h_3,h_4,f$ and $c$, we look for functions $h_1,h_2$ such that once inserted into the equation \eqref{oe1}, it yields that there exists a unique solution of the corresponding equation \eqref{oe1}  which satisfies the exponential asymptotic decay in the $L^2(0,\pi)$-norm \eqref{s1}.  Concerning the stability of equations like \eqref{oe1},  in  \cite{c17}, the authors provide some sufficient conditions on the parameters which assure the exponential decay of the solution. Besides this, in \cite{c18} the authors prove that for some very special cases of boundary disturbances, the backstepping stabilizing control still assures the stability of the system. The present paper provides an adaptive stabilizing control for equation \eqref{oe1}. To the best of our knowledge,   it represents a first result for the problem of stabilization of 1-D parabolic equations with  nonlocal boundary conditions as those in \eqref{oe1}. Anyway, there are some results of stabilization or controllability, but they are associated to hyperbolic equations with nonlocal boundary conditions, see \cite{c10} and \cite{c20}. 

The method we shall use here relies on the results obtained in \cite{ion1}. There, proportional type feedback stabilizing actuators were constructed for parabolic type equations with regular boundary conditions. Here, we show, in fact, that the same kind of controllers assure the stability for parabolic equations with nonlocal boundary conditions as-well. Anyway the proof is not straightforward and it implies some important changes in the method described in \cite{ion1}. The very simple explicit form of the controllers in \cite{ion1} allowed obtaining other first results in the literature concerning the stabilization of stochastic PDEs, see \cite{ 
ion1A,ion8,ion112}, or concerning the Cahn-Hilliard equations, see \cite{colli}. Others are related to the Navier-Stokes equations and MHD equations in a channel, see \cite{ion7,ion6}, or the phase-field system, see \cite{ion3}. For further details see the book \cite{book}.

In order to clarify the procedure, we start with a particular case, then we shall discuss different other cases (see Section \ref{s}). Let
$$h_1(x)=\int_0^xw(\xi)d\xi-H(x) \text{ and }h_2\equiv 0,$$
$$h_3(x)=\alpha H(x-\pi) \text{ and }h_4=H(x-\pi)+H(x).$$ Here, $H$ is the Heaviside function and $\alpha>0$ is some constant.
In this case, \eqref{oe1} reads as
\begin{equation}\label{o1} \left\lbrace \begin{array}{l}y_t(t,x)-y''(t,x)-cy(t,x)=f(t,x,y(t,x)),\ t>0, x\in(0,\pi),\\
\\
y(t,0)=\int_0^\pi w(x)y(t,x)dx,\ y'(t,0)+y'(t,\pi)+\alpha y(t,\pi)=0,\ t>0,\\
\\
y(0,x)=y_o(x),\ x\in(0,\pi).
\end{array} \right. \end{equation}Here, $w$ is seen as a control.

\section{Stabilization of the equation \eqref{o1}}

In the following, we denote by $L^2(0,\pi)$ the set of all square Lebesgue integrable functions on $(0,\pi)$. We set $\|\cdot\|$ and $\left<\cdot,\cdot\right>$ for the classical norm and scalar product in $L^2(0,\pi)$, respectively. Also, we denote by $\left<\cdot,\cdot\right>_{M}$ the euclidean scalar product in $\mathbb{R}^M$. For $k\in\mathbb{N}$, $H^k(0,\pi)$ will stand for the classical Sobolev space.

In a classical manner, we start with the stabilization  of the linearized system, given by
\begin{equation}\label{e1} \left\lbrace \begin{array}{l}y_t(t,x)-y^{\prime\prime}(t,x)-cy(t,x)=0,\ t>0,  x\in(0,\pi),\\
\\
y(t,0)=\int_0^\pi w(x)y(t,x)dx,\ y^\prime(t,0)+y^\prime(t,\pi)+\alpha y(t,\pi)=0,\\
\\
y(0,x)=y_o(x),\ x\in(0,\pi).
\end{array} \right. \end{equation} 

As mentioned above, we shall apply the results in \cite{ion1}, i.e., we shall look for a feedback controller $u=u(y)$ to stabilize the equation
\begin{equation}\label{ie1} \left\lbrace \begin{array}{l}y_t(t,x)-y^{\prime\prime}(t,x)-cy(t,x)=0,\ t>0,  x\in(0,\pi),\\
\\
y(t,0)=u(y(t)),\ y^\prime(t,0)+y^\prime(t,\pi)+\alpha y(t,\pi)=0,\ t>0,\\
\\
y(0,x)=y_o(x),\ x\in(0,\pi).
\end{array} \right. \end{equation}Then, since the controller $u$ will be given in a special proportional form, we shall recover the control function $w$, see \eqref{mio27} below.

We set $\mathbb{A}:\mathcal{D}(\mathbb{A})\subset L^2(0,\pi)\rightarrow L^2(0,\pi),$
$$\mathbb{A}y:=-y''-cy,\ \forall y\in\mathcal{D}(\mathbb{A}),$$where
$$\mathcal{D}(\mathbb{A})=\left\lbrace y\in H^2(0,\pi):\ y(0)=0,\ y'(0)+y'(\pi)+\alpha y(\pi)=0\right\rbrace.$$  By \cite{1}, we know that $\mathbb{A}$ has a countable set of eigenvalues $\left\lbrace \lambda_j\right\rbrace_{j=0}^\infty $ described as follows
$$\lambda_j=\left\lbrace \begin{array}{ll}(2k+1)^2-c & ,\text{if } j=2k, \ k\in\mathbb{N},\\ (2\beta_k)^2-c & ,\text{if } j=2k+1,\ k\in\mathbb{N} .\end{array} \right. $$ Here, $\beta_k,\ k\in\mathbb{N},$ are the roots of the equation 
$$\cot(\beta\pi)=-\frac{\alpha}{2\beta},$$ for which we know that they satisfy  
$$2k+1 < 2\beta_k< 2k+2,\ k=0,1,2,3,....$$ Easily seen, given $\rho>0$, there exists $N\in\mathbb{N}$ such that 
$$\lambda_0<\lambda_1<...<\lambda_{2N+1}\leq\rho<\lambda_{2N+2}<\lambda_{2N+3}<....$$The corresponding eigenfunctions are precisely given in \cite{1}. More precisely they are $\left\lbrace w_{k1},\ w_{k2}\right\rbrace_{k=0}^\infty$, where
$$w_{k1}=\sin((2k+1)x),\  \text{and }w_{k2}=\sin(2\beta_kx), \ k\in\mathbb{N}.$$ As stated and proved in \cite[Lemma 4.1]{1}, it happens that the above system does not form a basis in $L^2(0,\pi)$.  That is why, in \cite{1}, the authors introduced the following new set of functions
\begin{equation}\label{e3}\varphi_j(x)=\left\lbrace \begin{array}{ll} w_{k1}(x)&,\text{if }j=2k,\\
\left[ w_{k2}(x)-w_{k1}(x)\right](2\delta_k)^{-1}&,\text{if }j=2k+1,\end{array} \ \ \ \  k\in\mathbb{N}.\right.  \end{equation}Here, $\delta_k:=\beta_k-k-\frac{1}{2}.$ Then, in \cite[Lemma 5.1]{1} they proved that the system $\left\lbrace \varphi_j\right\rbrace _{j=0}^\infty$ forms a Riesz basis in $L^2(0,\pi)$. Moreover, they do also precised the bi-orthonormal system to $\left\lbrace \varphi_j\right\rbrace _{j=0}^\infty$, which is given by
\begin{equation}\label{e4}\psi_j(x)=\left\lbrace \begin{array}{ll}v_{k2}(x)+v_{k1}(x)&,\text{if }j=2k,\\
2\delta_kv_{k2}(x)&,\text{if }j=2k+1,\end{array}  \ \ \ \ \ k\in\mathbb{N}. \right. \end{equation} Where
$$v_{k1}(x)=\frac{2}{\pi}\left\lbrace \sin((2k+1)x)-\frac{2k+1}{\alpha}\cos((2k+1)x)\right\rbrace ,\ k=0,1,2,...,$$and
$$v_{k2}=C_{k2}\left\lbrace \sin(2\beta_kx)-\frac{2\beta_k}{\alpha}\cos(2\beta_kx)\right\rbrace,\ k=0,1,2,...$$Here, $C_{2k}$ is some constant which assures that the systems $ \left\lbrace w_{k1},w_{k2}\right\rbrace_{k=0}^\infty $ and $\left\lbrace v_{k1},v_{k2}\right\rbrace_{k=0}^\infty$  are bi-orthonormal. 

It is easy to see that, at the boundary, we have
$$\psi_j0)+\psi_j(\pi)=0 \text{ and } \psi_j'(\pi)+\alpha\psi_j(\pi)=0,\ j\in\mathbb{N}.$$

In order to lift the boundary controller into the equations, let us introduce the following map, which is usually called the Dirichlet map. For $\gamma>0$, let $D$ be the solution to the equation 
\begin{equation}\label{e5}\left\lbrace \begin{array}{l}\begin{aligned}&-D''(x)-cD(x)-2\sum_{j=0}^{2N+1}\lambda_j\left<D,\psi_j\right>\varphi_j\\&
-2\sum_{j=0}^{N}(2\beta_j+2j+1)\left<D,\psi_{2j+1}\right>\varphi_{2j}+\gamma D=0,\ x\in(0,\pi),\end{aligned}\\
\\
D(0)=1,\ D'(0)+D'(\pi)+\alpha D(\pi)=0.\end{array} \right. \end{equation} In the lemma below, we show that, for $\gamma$ large enough, the above equation is well-posed.
\begin{lemma}
	If $\gamma>0$ is large enough, then equation \eqref{e5} has a unique solution.
\end{lemma}
\begin{remark} Notice that, in the present case, the lifting operator $D$ is defined in a different way then that one used in \cite[Eq. (2.1)]{ion1}. The reason is because the system $\left\lbrace \varphi_j\right\rbrace_{j=0}^{\infty}$ is not an eigen-basis of $L^2(0,\pi)$. Thus, after the projection,   the diagonal matrix  $\Lambda$ in  \cite[Eq.(4.11)]{ion1}, has, in the present case, the non-diagonal form \eqref{e30} below.
	\end{remark} 
\begin{proof}We shall rely on the results obtained in \cite{2}. More exactly, as in the proof of \cite[Lemma 3.1]{2}, we begin with considering the second order differential equation
	\begin{equation}\label{e7}\left\lbrace \begin{array}{l}-v''(x)+(\gamma-c)v(x)=\chi(x),\ x\in(0,\pi),\\
	v(0)=0, \ v(\pi)=\omega. \end{array}\right. \end{equation} 
	If $\gamma>c$ and $\chi$ is twice differentiable, then if equation \eqref{e7} has a solution it can be represented in an integral form  as
$$\begin{aligned} &v(x)=\frac{1}{1-e^{-2\pi b}}\left(  e^{-(\pi-x)b}-e^{-(\pi+x)b}\right) \omega   \\&
	-\frac{1}{1-e^{-2\pi b}}\left\lbrace \frac{1}{2b}(e^{-(\pi-x)b}-e^{-(\pi+x)b})\int_0^\pi (e^{-(\pi-\xi)b}-e^{-(\pi+\xi)b})\chi(\xi)d\xi\right\rbrace\\&
	+\frac{1}{2b}\int_0^\pi\left( e^{-|x-\xi|b}-e^{-(x+\xi)b}\right)\chi(\xi)d\xi.  \end{aligned}$$Here, $b=\sqrt{\gamma-c}.$
		
Imposing that $v'(0)+v'(\pi)+\alpha v(\pi)+\alpha=0$ and arguing as in \cite[Eqs. (13)-(14)]{1}, we deduce that the solution $v$ can be expressed as
$$v(x)=\int_0^\pi G(x,\xi)\chi(\xi)d\xi,$$where
$$\begin{aligned}G(x,\xi)=& (1-e^{-2\pi b})^{-1}\left(e^{-(\pi-x)b}-e^{-(\pi+x)b}+\alpha \right) )P\left(e^{-(\pi-\xi)b}+e^{-\xi b} \right)\\&
+(1-e^{-2\pi b})^{-1}\left( e^{-(\pi-x)b}-e^{-(\pi+x)b}) \right) \frac{1}{2b}\left( e^{-(\pi-\xi)b}+e^{-(\pi+\xi)b}\right)\\&
+\frac{1}{2b}(1-e^{-2\pi b})^{-1}\left(1-e^{-2\pi b}\right)\left( e^{-|x-\xi|b}-e^{-(x+\xi)b}\right),   \end{aligned}$$where
$$P=\left( \alpha +b\frac{(1+e^{-\pi b})^2}{1-e^{-2\pi b}}\right) ^{-1}.$$

Now, let us replace $\chi$ by
$$\begin{aligned}\chi=&2\sum_{j=0}^{2N+1}\lambda_j\left<v,\psi_j\right>\varphi_j+
2\sum_{j=0}^{N}(2\beta_j+2j+1)\left<v,\psi_{2j+1}\right>\varphi_{2j}\\&
+2\sum_{j=0}^{2N+1}\lambda_j\left<1,\psi_j\right>\varphi_j+
2\sum_{j=0}^{N}(2\beta_j+2j+1)\left<1,\psi_{2j+1}\right>\varphi_{2j}-(\gamma-c).\end{aligned}$$ Hence, one may define the operator $T:L^2(0,\pi)
\rightarrow L^2(0,\pi)$ as 
$$\begin{aligned}&T(v)(x):=\\&
2\int_0^\pi G(x,\xi)\left[\sum_{j=0}^{2N+1}\lambda_j\left<v,\psi_j\right>\varphi_j(\xi)+
\sum_{j=0}^{N}(2\beta_j+2j+1)\left<v,\psi_{2j+1}\right>\varphi_{2j}(\xi)\right]d\xi\\&
+2\int_0^\pi G(x,\xi)\left[ 2\sum_{j=0}^{2N+1}\lambda_j\left<1,\psi_j\right>\varphi_j(\xi)+
2\sum_{j=0}^{N}(2\beta_j+2j+1)\left<1,\psi_{2j+1}\right>\varphi_{2j}(\xi)\right] d\xi\\&
-(\gamma-c)\int_0^\pi G(x,\xi) d\xi. \end{aligned}$$Thus, the solution $v$ to \eqref{e7} for the corresponding $\chi$, if  exists, it is a fixed point of $T$. 

Taking into account that all the exponentials, appearing in the definition of $G$, have negative power, due to the definition of $P$ and $b$, we get that
$$|G(x,\xi)|\leq C\frac{1}{\sqrt{\gamma}},\ \forall x,\xi \in (0,\pi),$$for some postive constant $C$, independent of $\gamma$. Moreover, since all the functions $\left\lbrace \varphi_j\right\rbrace_{j=0}^{2N+1}, \ \left\lbrace \psi_j\right\rbrace_{j=0}^{2N+1} $ are bounded, we get via the Schwarz inequality that
$$\left| \sum_{j=0}^{2N+1}\lambda_j\left<v,\psi_j\right>\varphi_j(\xi)+
\sum_{j=0}^{N}(2\beta_j+1)\left<v,\psi_{2j+1}\right>\varphi_{2j}(\xi)\right| \leq C\|v\|,$$ $\forall \xi\in (0,\pi)$, for some positive constant $C$, independent of $\gamma$. 

 Let any $v_1,v_2\in L^2(0,\pi)$, it follows from the above estimates that
$$\|Tv_1-Tv_2\|\leq C\frac{1}{\sqrt{\gamma}}\|v_1-v_2\|,$$ for some positive constant $C$ independent of $\gamma$ and $v_1,v_2$. Thus, if $\gamma$ is large enough, it yields that $T$ is a contraction. Consequently, via the contraction mapping theorem it follows that $T$ has a unique fixed point, which is the solution $v$ for the corresponding $\chi$. To conclude, setting $D:=v+1$, we immediately get that $D$ is the unique solution to \eqref{e5}.

\end{proof}

For latter purpose, let us perform the computations for 
$$\left<D,\psi_k\right>,\ k=0,1,2,...,2N+1,$$where $D$ is the unique solution of \eqref{e5}. To this end, taking into account the bi-orthogonality of the systems $\left\lbrace \varphi_j\right\rbrace_{j=0}^\infty $ and  $\left\lbrace \psi_j\right\rbrace_{j=0}^\infty $ together with their boundary conditions,  it yields, after scalar multiplication of equation \eqref{e5} by $\psi_{2j}$ and integration by parts, that
\begin{equation}\label{e10}(\gamma-\lambda_{2j})\left<D,\psi_{2j}\right>-(2\beta_j+2j+1)\left<D,\psi_{2j+1}\right>=\psi'_{2j}(0),\end{equation}since $-\psi''_{2j}-c\psi_{2j}=\lambda_{2j}\psi_{2j}+(2\beta_j+2j+1)\psi_{2j+1}$. On the other hand, if we scalarly multiply equation \eqref{e5} by $\psi_{2j+1}$, we get that
 \begin{equation}\label{e11} \left<D,\psi_{2j+1}\right>=\frac{1}{\gamma-\lambda_{2j+1}}\psi'_{2j+1}(0),\end{equation}since $-\psi''_{2j+1}-c\psi_{2j+1}=\lambda_{2j+1}\psi_{2j+1}.$ 
 
 In virtue of \eqref{e10} and \eqref{e11}, we conclude  that
 \begin{equation}\label{e15}\left\lbrace \begin{array}{l}\left<D,\psi_{2j}\right>=\frac{1}{\gamma-\lambda_{2j}}l_{2j}:= \frac{1}{\gamma-\lambda_{2j}}\left[\psi'_{2j}(0)+\frac{2\beta_j+2j+1}{\gamma-\lambda_{2j+1}}\psi'_{2j+1}(0) \right],\\
 \\
 \left<D,\psi_{2j+1}\right>=\frac{1}{\gamma-\lambda_{2j+1}}l_{2j+1}:=\frac{1}{\gamma-\lambda_{2j+1}} \psi'_{2j+1}(0),\end{array}\right.  \end{equation}for $j=0,1,2,...,N$. It is easy to see that since $\psi'_{2j+1}(0)\neq0$,  for $\gamma$ properly chosen, we have
 $$l_{k}\neq0, \ k=0,1,,...2N+1.$$ 
 
 Let $\rho<\gamma_0<\gamma_1<...<\gamma_{2N+1}$ be $2N+2$ large enough positive numbers such that, for each of them the equation \eqref{e5} is well-posed, and denote by $D_{\gamma_k},\ k=0,1,...,2N+1,$ each corresponding solution.
 
 In the next lines we describe the boundary feedback controller which we claim that assures the exponential stability of the equation.
 We denote by $\textbf{B}$ the next square matrix 
 \begin{equation}\label{ie5}\textbf{B}:=\left( l_il_j\right)_{0\leq i,j\leq 2N+1} \end{equation}and multiply it on both sides by $$\Lambda_{\gamma_k}:=diag\left(\frac{1}{\gamma_k-\lambda_0},\frac{1}{\gamma_k-\lambda_1},...,\frac{1}{\gamma_k-\lambda_{2N+1}}\right), \ k=0,1,...,2N+1,$$ to define
 \begin{equation}\label{ie8}B_k:=\Lambda_{\gamma_k}\textbf{B}\Lambda_{\gamma_k},\ k=0,...,2N+1.\end{equation}Then, we introduce the matrix
 \begin{equation}\label{ie9}A:=(B_0+B_1+...+B_{2N+1})^{-1}\end{equation}which, as claimed in \cite[Lemma 5.2]{ion1}, is well-defined since $l_k\neq0,\ k=0,1,...,2N+1$. Then, for each $k=0,1,...,2N+1$, we set the following feedback forms
 \begin{equation}\label{ie10}u_k(y):=\left<\Lambda_{\gamma_k}A\left(\begin{array}{c}\left<y,\psi_0\right>\\ \left<y,\psi_1\right>\\\ddots \\
 \left<y,\psi_{2N+1}\right>\end{array}\right), \left(\begin{array}{c}l_0\\ l_1\\
 \ddots \\
 l_{2N+1}
 \end{array}\right)\right>_{2N+2}, \end{equation} then, introduce $u$ as 
 \begin{equation}\label{mio27}\begin{aligned}u(y):=&-\left[ u_0(y)+u_1(y)+...+u_{2N+1}(y)\right] \\&
 \\&
 =-\left<\sum_{k=0}^{2N+1}\Lambda_{\gamma_k}A\left(\begin{array}{c}\left<y,\psi_0\right>\\ \left<y,\psi_1\right>\\\ddots \\
 \left<y,\psi_{2N+1}\right>\end{array}\right), \left(\begin{array}{c}l_0\\ l_1\\
 \ddots \\
 l_{2N+1}
 \end{array}\right)\right>_{2N+2}.\end{aligned}\end{equation} It is easy to see that the control $u$ can be equivalently rewritten as
 \begin{equation}\label{e40}u(y)=-\int_0^\pi \left<\sum_{k=0}^{2N+1}\Lambda_{\gamma_k}A\left(\begin{array}{c}\psi_0(x)\\ \psi_1(x)\\
 \ddots \\
 \psi_{2N+1}(x)
 \end{array}\right),\left(\begin{array}{c}l_0\\ l_1\\
 \ddots \\
 l_{2N+1}
 \end{array}\right)\right>_{2N+2} \times y(x)dx.\end{equation}
 
Let us return to equation \eqref{ie1}, in which we plug the controller $u$, given by \eqref{mio27}, i.e.
\begin{equation}\label{iie1} \left\lbrace \begin{array}{l}y_t(t,x)-y''(t,x)-cy(t,x)=0,\ x\in(0,\pi),\\
\\
y(t,0)=-\left[ u_0(y(t))+u_1(y(t))+...+u_{2N+1}(y(t))\right] ,\\
\\
\ y'(t,0)+y'(t,\pi)+\alpha y(t,\pi)=0,\ t>0,\\
\\
y(0,x)=y_o(x),\ x\in(0,\pi).
\end{array} \right. \end{equation}Next define the new variable
$$z(t,x):=y(t,x)+\sum_{k=0}^{2N+1}u_k(t)D_{\gamma_k}(x),\ t\geq0,\ x\in(0,\pi).$$ We get, likewise in \cite[Eqs.(4.7)-(4.9)]{ion1}, that  equation \eqref{iie1} can be rewritten in terms of $z$ as
 \begin{equation}\label{e24}\left\lbrace \begin{array}{l}z_t+\mathbb{A}z=\mathcal{R}(\left<z,\psi_0\right>,...,\left<z,\psi_{2N+1}\right>), \ t>0,\\
 \\
 z(0)=z_o:=y_o+\sum_{k=0}^{2N+1}u_k(0)D_{\gamma_k}.\end{array}\right. \end{equation} Here,
 $$\begin{aligned}&\mathcal{R} (\left<z,\psi_0\right>,...,\left<z,\psi_{2N+1}\right>) :=\\&
 =-\left( \sum_{k=0}^{2N+1}u_k(t)D_{\gamma_k}(x)\right)_t-2\sum_{k,j=0}^{2N+1}\lambda_k\left<u_k(t)D_{\gamma_k},\psi_j\right>\varphi_j(x)\\&
 -2\sum_{k=0,j=0}^{k=2N+1,j=N}(2\beta_j+2j+1)\left<u_k(t)D_{\gamma_k},\psi_{2j+1}\right>\varphi_{2j}(x)+\sum_{k=0}^{2N+1}\gamma_ku_k(t)D_{\gamma_k}(x). 
 \end{aligned}$$
 Concerning the stability of the linear equation  we have the following result.
 \begin{theorem}\label{t1}The equation \eqref{e24} is well-posed, and its unique solution satisfies
 \begin{equation}\label{e27}\|z(t)\|\leq Ce^{-\rho t}\|z_o\|,\ \forall t\geq0,\end{equation}where $C$ is some positive constant. 

Consequently, by \eqref{e40}, once we plugg the function
$$\omega(x)=- \left<\sum_{k=0}^{2N+1}\Lambda_{\gamma_k}A\left(\begin{array}{c}\psi_0(x)\\ \psi_1(x)\\
\ddots \\
\psi_{2N+1}(x)
\end{array}\right),\left(\begin{array}{c}l_0\\ l_1\\
\ddots \\
l_{2N+1}
\end{array}\right)\right>_{2N+2} $$ into the linear equation \eqref{e1},  it yileds that the unique corresponding solution of the closed-loop equation \eqref{e1} satisfies
\begin{equation}\label{e41}\|y(t)\|\leq Ce^{-\rho t}\|y_o\|,\ \forall t\geq0,\end{equation}for some positive constant $C$.\end{theorem}
\begin{remark}Note that the exponential decay \eqref{e41} is arbitrarly fast since $\rho>0$ was arbitrarily chosen. Anyway, the faster the decay is, the bigger the $N$ is, and consequently, the bigger the dimension of the controller is.
	\end{remark}
\begin{proof} Concerning the well-posedness of equation \eqref{e24}, which is with nonlocal boundary conditions, one may argue as in \cite{2}.
	
	 Next, since $\left\lbrace \varphi_j\right\rbrace_{j=0}^\infty $ is a Riesz basis in $L^2(0,\pi)$, with $\left\lbrace \psi_j\right\rbrace_{j=0}^\infty $ its bi-orthonormal system, we have that the solution $z$ can be decomposed as
	 $$z=z_N+\zeta_N,$$where
	 $$z_N:=\sum_{j=0}^{2N+1}\left<z,\psi_j\right>\varphi_j \text{ and }\zeta_N:=\sum_{j=2N+2}^\infty\left<z,\psi_j\right>\varphi_j.$$ Next we shall apply the classical decomposition method. It is not a spectral decomposition since $\varphi_j$ are not eigenfunctions of $\mathbb{A}$, but combinations of them.  Recall the form of $\varphi_j$'s given by \eqref{e3}. Then, simple computations, lead to
	 $$\mathbb{A}\varphi_j=\left\lbrace \begin{array}{ll} \lambda_{2k}\varphi_{2k}&,\text{if }j=2k,\\
	 \lambda_{2k+1}\varphi_{2k+1}+(2\beta_k+2k+1)\varphi_{2k}&,\text{if }j=2k+1,\end{array}\  \ \ \ \ k=0,1,2,.....\right.  $$ Hence, using the bi-orthonormality of the systems $\left\lbrace \varphi_j\right\rbrace_{j=0}^\infty $ and $\left\lbrace \psi_j\right\rbrace_{j=0}^\infty $, we obtain that  $$\left(\begin{array}{c} \left<\mathbb{A}z,\psi_0\right>\\ \left<\mathbb{A}z,\psi_1\right>\\ \ddots \\ \left<\mathbb{A}z,\psi_{2N+1}\right>\end{array}\right)=\Lambda\mathcal{Z},$$where
	 \begin{equation}\label{e30}\Lambda:=
	 \left(\begin{array}{ccccccc}\lambda_0&2\beta_0+1&0&\dots&0&0&0\\
	 0&\lambda_1&0&\dots&0&0&0\\
	 0&0&\lambda_2&2\beta_1+3&0&\dots&0\\
	 \ddots&\ddots&\ddots&\ddots&\ddots&\ddots&\ddots\\
	 0&0&0&\dots&\lambda_{2N-1}&0&0\\
	 0&0&0&0&\dots&\lambda_{2N}&2\beta_N+2N+1\\
	 0&0&0&\dots&0&0&\lambda_{2N+1} \end{array}\right).\end{equation}Above, 
	 $$\mathcal{Z}=\left( \begin{array}{c}\left<z,\psi_0\right>\\ \left<z,\psi_1\right>\\ \ddots \\ \left<z,\psi_{2N+1}\right> \end{array}\right).$$
	 
	 Similar arguments as in \cite[Eqs.(4.7)-(4.13)]{ion1} lead us to the fact that the projection of the equation \eqref{e24} on the space spaned by $\left\lbrace\varphi_j \right\rbrace_{j=0}^{2N+1}$ is given by 
 $$\mathcal{Z}_t+\Lambda\mathcal{Z}=\frac{1}{2}\mathcal{Z}+\Lambda\mathcal{Z}-\frac{1}{2}\sum_{k=0}^{2N+1}\gamma_kB_kA\mathcal{Z},$$ or, equivalently
  $$\begin{aligned}\mathcal{Z}_t&=-\sum_{k=0}^{2N+1}\gamma_kB_kA\mathcal{Z}\\&
  =-\gamma_0\sum_{k=0}^{2N+1}B_kA \mathcal{Z}+\sum_{k=1}^{2N+1}(\gamma_0-\gamma_k)B_k A\mathcal{Z}\\&
  =-\gamma_0\mathcal{Z}+\sum_{k=1}^{2N+1}(\gamma_0-\gamma_k)B_k A\mathcal{Z},\ t>0,\end{aligned}$$ since $A$ is the inverse of the sum of $B_k'$s. Then, arguing as in \cite[Eq.(4.13)]{ion1}, since $\left<A\cdot,\cdot\right>_{2N+2}$ is a Lyapunov function for the the above equation, we deduce that $\mathcal{Z}$ satisfies
  $$\|\mathcal{Z}(t)\|_{\mathbb{R}^{2N+2}}\leq Ce^{-\gamma_0t}\|\mathcal{Z}(0)\|_{\mathbb{R}^{2N+2}},\ \forall t\geq0.$$ From where it follows the $L^2-$exponential decay of the component $z_N$, i.e.,
  $$\|z_N(t)\|\leq Ce^{-\rho t}\|z_N(0)\|,\ \forall t\geq0,$$since $\rho<\gamma_0$ and $\left\lbrace \varphi_j\right\rbrace_{j=0}^\infty $ is a Riesz basis in $L^2(0,\pi)$. Finally, the component $\zeta_N$ is also exponentially asymptotically decaying in the $L^2-$norm since the projection of the operator $\mathbb{A}$, on the space spanned by $\left\lbrace \varphi_j\right\rbrace_{j=2N+2}^\infty $, involves  the eigenvalues $\lambda_j>\rho,\ j=2N+2,2N+3,...,$. Thus, if $\rho$ is large enough (that is, $N$ is large enough), one can use similars ideas as in \cite{c17} to show the asymptotic exponential stability of the equation corresponding for $\zeta_N$.  The rest of the proof follows identically as the proof of \cite[Theorem 2.1]{ion1}. The details are omitted.
 
 \end{proof}

 Finally, arguing as in the proof of \cite[Theorem 2.2]{ion1}, we deduce the following result concerning the stabilization of the nonlinear system \eqref{o1}.
 \begin{theorem}\label{t2} Once we plugg the function
 	$$\omega(x)=-\left<\sum_{k=0}^{2N+1}\Lambda_{\gamma_k}A\left(\begin{array}{c}\psi_0(x)\\ \psi_1(x)\\
 	\ddots \\
 	\psi_{2N+1}(x)
 	\end{array}\right),\left(\begin{array}{c}l_0\\ l_1\\
 	\ddots \\
 	l_{2N+1}
 	\end{array}\right)\right>_{2N+2} $$into the equation \eqref{o1}, it yields the existence of some $\delta>0$ such that, for each initial data $\|y_o\|\leq \delta$, there exists a unique solution to the corresponding closed-loop equation \eqref{o1}, which satisfies
 	$$\|y(t)\|\leq Ce^{-\rho t}\|y_o\|,\ \forall t\geq0,$$ for some positive constant $C$. 
 \end{theorem}

\section{Other examples of boundary nonlocal conditions}\label{s}
Let us cover some other cases
\begin{itemize}
	 \item[1)] In equation \eqref{oe1}, let us consider
$$h_1(x)=\int_0^x\omega(\xi)d\xi-H(x); h_2\equiv 0;$$
$$\ h_3(x)=\alpha H(x-\pi);\ h_4(x)=H(x-\pi)-H(x),$$with $\alpha\in\mathbb{R}$.This gives the B.C.
$$y(t,0)=\int_0^\pi\omega(x)y(t,x)dx\text{ and }y'(t,0)=y'(t,\pi)+\alpha y(t,\pi).$$ In \cite{c1} the author describes the eigenfunctions of the corresponding linear governing operator $\mathbb{A}$, with the B.C.
$$w(0)=0,\ w'(0)=w'(\pi)+\alpha w(\pi).$$ They are similar with the ones above. So, arguing as in \cite{1}, one may construct from the eigenfunctions system a Riesz basis of $L^2(0,\pi)$. Then define the Dirichlet lifting operator $D$ as in \eqref{e5}, and compute the products $\left<D,\psi_j\right>,\ j=0,1,2,...,2N+1$, introducing this way the $l_j$'s. Again, by proper choice of $\gamma'$s, one obtains that  $l_j\neq0,\ \forall j$. Therefore,  a similar stabilizer as that one in \eqref{e40} can be constructed for the above nonlocal B.C.

\item[2)] In case
$$h_1(x)=\int_0^x\omega(\xi)d\xi-H(x);\  h_2\equiv0;$$
$$h_3(x)=\alpha x- H(x-\pi);\ h_4\equiv0.$$ This gives the B.C.
$$y(t,0)=\int_0^\pi\omega(x)y(t,x)dx \text{ and }y(\pi)=\alpha\int_0^\pi y(t,x)dx.$$The corresponding linear operator $\mathbb{A}$ has the B.C.
$$w(0)=0,\ w(\pi)=\alpha \int_0^\pi w(x)dx.$$Computing the eigenfunctions, we get
$$\sin(\sqrt{\lambda}\pi)=\alpha\int_0^\pi \sin(\sqrt{\lambda}x)dx,$$which leads to two families of eigenfunctions
$$w_{1k}(x)= \sin(2kx) \text{ and }w_{2k}= \sin(\beta_kx),$$where $\beta_k$ is the root of
$$\cot(\beta\pi)=\frac{\alpha}{2\beta}.$$ But these are again similar with the eigenfunctions from the above section (see also pp. 983-984 in \cite{1}). Therefore, the same procedure of the construction of the Riesz basis in \cite{1} can be performed and a similar feedback stabilizer can be constructed in this case as-well.
\item[3)] Let us consider now the case 
$$h_1(x)=\int_0^x\omega(\xi)d\xi-H(x);\  h_2\equiv0;$$
$$h_3(x)=H(x-\eta);\ h_4\equiv 0,$$where $\eta\in(0,\pi).$ The B.C. read as
$$y(t,0)=\int_0^\pi\omega(x)y(t,x)dx \text{ and }y(\eta)=0.$$ In \cite{sayed}, the authors proved that the corresponding eigenfunctions of the linear operator form an orthonormal eigenbasis of $L^2(0,\pi)$, and the eigenfunction expansion formula holds true. Consequently, the same stabilizing control design technique can be applied, directly on the eigenfunction system.
\end{itemize}
In conclusion of this section, we emphasize that there can be found many other examples of semilinear heat equations with nonlocal boundary conditions, to which the above stabilizing control design method works equally-well. 
\section{Conclusions} Based on the results in \cite{ion1}, we designed here boundary stabilizing controllers for the 1-D semilinear heat equation with nonlocal boundary conditions. To clearify the method, we considered some particular cases of B.C. only. Anyway, the method works  for many other cases, provided that  given the Laplace operator with nonlocal B.C.,  one  solves the associated spectral problem, then constructs from the eigenfunction system a Riesz base in $L^2$ and a bi-orthonormal system, for which, the first $M$  functions satisfy some unique continuation property derived from the scalar product between them and the lifting Dirichlet operator. The form of the actuator involves only the first $M$ eigenfunctions of the linear operator. Regarding this, we stress that in \cite{numeric}, the author provides an efficent method to compute the spectrum of the Laplace operator with nonlocal boundary conditions of the type we used throughout this paper.

\end{document}